\documentclass[reqno]{amsart}
\usepackage{amsmath}

\newcommand{\RR}{\ensuremath{\mathbb{R}}}
\newcommand{\HH}{\ensuremath{\mathbb{H}}}

\newcommand{\hgot}{\ensuremath{\mathfrak{h}}}

\newcommand{\ggot}{\ensuremath{\mathfrak{g}}}
\newcommand{\mgot}{\ensuremath{\mathfrak{m}}}
\newcommand{\xgot}{\ensuremath{\mathfrak{X}}}

\DeclareMathOperator{\Ker}{Ker} 
\DeclareMathOperator{\ident}{id}
\DeclareMathOperator{\Image}{Im} 
\DeclareMathOperator{\Tr}{Tr}
\DeclareMathOperator{\real}{Re} 
\DeclareMathOperator{\md}{mod}

\newtheorem{thm}{Theorem}
\newtheorem{as}{Assumption}

\newtheorem{lem}{Lemma}
\newtheorem{prop}{Proposition}

\theoremstyle{definition}
\newtheorem{defn}{Definition}

\begin{document}

\title[Invariant metric $f$-structures]{Invariant metric $f$-structures on specific homogeneous reductive spaces}

\author[A. Sakovich]{Anna Sakovich}

\address{Anna Sakovich\newline Faculty of Mathematics and Mechanics\newline
Belarusian State University\newline Nezavisimosti av.~4\newline Minsk
220050, BELARUS }


\email{anya\_sakovich@tut.by}

\thanks{{\it Acknowledgements.} The author is grateful to Vitaly V. Balashchenko for stimulating discussion that has initiated this paper.}






\maketitle

\begin{abstract}
For homogeneous reductive spaces $G/H$ with reductive complements decomposable into an orthogonal sum $\mathfrak{m}=\mathfrak{m} _1 \oplus \mathfrak{m} _2 \oplus \mathfrak{m} _3$  of three $Ad(H)$-invariant irreducible mutually inequivalent submodules we establish simple conditions under which an invariant metric $f$-structure $(f,g)$ belongs to the classes ${\bf G_1 f}$, ${\bf NKf}$, and ${\bf Kill \, f}$ of generalized Hermitian geometry. The statements obtained are then illustrated with four examples. Namely we consider invariant metric $f$-structures on the manifolds of oriented flags $SO(n)/SO(2)\times SO(n-3)$ $(n \geq 4)$, the Stiefel manifold $SO(4)/SO(2)$, the complex flag manifold $SU(3)/T_{max}$, and the quaternionic flag manifold $Sp(3)/SU(2)\times SU(2) \times SU(2)$.
\end{abstract}

\section*{ Introduction}

The concept of generalized Hermitian geometry (see, for example, \cite{Ki1}) was created in the 1980s as a natural consequence of the development of Hermitian geometry and the theory of almost contact structures. One of the central objects in this concept is the metric $f$-structure $(f,g)$, that is, an $f$-structure \cite{Y} $f$ compatible with an invariant Riemannian metric $g$. 

An interesting problem that arises in this context is to determine whether a given metric $f$-structure belongs to the main classes of generalized Hermitian geometry, for example, to the classes ${\bf G_1 f}$ (see \cite{Ki1}), ${\bf NKf}$ (see \cite{B3} and \cite{B6}),  and ${\bf Kill \, f}$ (see \cite{Gr1} and \cite{Gr2}). It should be emphasized that in the case of naturally reductive manifolds \cite{KN} there exist a  number of results  that transform this problem into an easy computational task (\cite{B3}, \cite{B2}, \cite{B4}, \cite{B5}). However, in the case of an arbitrary Riemannian metric this problem is not an easy one, at least because it involves the calculation of the implicitly defined Levi-Civita connection.

In this paper we consider invariant metric $f$-structures $(f,g)$ on specific homogeneous reductive spaces $G/H$, namely on homogeneous reductive spaces that satisfy the following set of conditions:

\begin{itemize}
\item[1)] $G$ is a compact semisimple Lie group (hence the Killing form $B$ of $G$ is negative definite).
\item[2)] The reductive complement $\mgot$ admits the decomposition
\begin{equation*}
\mgot=\mgot _1 \oplus \mgot _2 \oplus \mgot _3
\end{equation*}
into the direct sum of $Ad(H)$-invariant irreducible mutually non-equivalent submodules and
this decomposition is $B$-orthogonal. 
\item[3)]
\begin{equation*}
0\neq [\mgot _i, \mgot _{i+1}] \subset \mgot _{i+2} \; (\md \;
3),\; i=1,2,3.
\end{equation*}
\item[4)]
\begin{equation*}
[\mgot _i, \mgot _i]\subset \hgot, \;i=1,2,3,
\end{equation*}
where $\hgot$ is the Lie algebra of $H$.
\end{itemize}

In this case it is not difficult to obtain an explicit formula for the Levi-Civita connection of a Riemannian manifold $(G/H,g)$. At the same time, for any nontrivial invariant $f$-structure which is not an almost complex structure \cite{KN} there exists such $i \in \{1,2,3\}$ that either $\Image f=\mgot _i$ or $\Ker f=\mgot _i$. This, in its turn, has enabled us to obtain easy-to-check characteristic conditions (Theorem 2 and Theorem 3) for metric $f$-structures $(f,g)$ under which they belong to the aforementioned classes of generalized Hermitian geometry.

Note that this paper was initiated by the study of the manifolds of oriented flags $SO(n)/SO(2)\times SO(n-3)$ $(n \geq 4)$. In \cite{IJMMS} it was shown that these homogeneous spaces satisfy the conditions 1) -- 4). In the last section of this paper we provide other examples of such spaces. Namely, by making use of Theorem 2 and Theorem 3, we consider invariant metric $f$-structures on the Stiefel manifold $SO(4)/SO(2)$, the complex flag manifold $SU(3)/T_{max}$, and the quaternionic flag manifold $Sp(3)/SU(2)\times SU(2) \times SU(2)$.

 \section{Preliminaries}
 
 \subsection{Invariant $f$-structures on homogeneous reductive spaces}

Homogeneous reductive spaces make up the main subject of our
further considerations. Therefore we begin with recollecting some basic facts related to them.

\begin{defn}\cite{N}
Let $G$ be a connected Lie group, $H$ its closed subgroup, $\ggot$
and $\hgot$ the corresponding Lie algebras. $G/H$ is called a {\it homogeneous 
reductive space} if there exists $\mgot \subset \ggot$
such that
\begin{itemize}
\item [1)] $\ggot=\hgot \oplus \mgot$.
\item [2)] For any $h \in H \; Ad(h) \mgot \subset \mgot$.
\end{itemize}
$\ggot=\hgot \oplus \mgot$ is the {\it reductive decomposition}
corresponding to $G/H$ and $\mgot$ is the {\it reductive
complement}.
\end{defn}

For any homogeneous reductive space $G/H$ its reductive complement $\mgot$ can be identified with the tangent space to $G/H$ at the
point $o=H$ in the following sense:
\begin{equation*}
\text{for any } h \in H \; d \tau (h)_o=Ad(h),\text{ where } \tau (g):G/H \rightarrow G/H, \;
xH \rightarrow (gx)H.
\end{equation*}

Since all homogeneous spaces to be discussed in this paper are
reductive, we agree to identify their reductive complements and their
tangent spaces at the point $o$.

An {\it affinor structure} on a smooth manifold is known to be a
tensor field of type (1,1) realized as a field of endomorphisms
acting on its tangent bundle. In this paper we will be primarily interested in the {\it almost complex structure} \cite{KN} (such an affinor structure $J$ that $J^2=-\ident$) and the {\it $f$-structure} \cite{Y} (an affinor structure $f$ satisfying $f^3+f=0$). 

\begin{defn}\cite{V} Let $G/H$ be a homogeneous manifold, $F$ an
affinor structure. $F$ is called {\it invariant} with respect to
$G$ if for any $g \in G$
\begin{equation*}
d \tau (g) \circ F = F \circ d \tau (g).
\end{equation*}
\end{defn}

It is known that any invariant affinor structure $F$ on a reductive homogeneous space $G/H$  is completely
determined by its value $F_o$ at the point $o=H$, where $F_o$ is a
linear operator on the reductive complement $\mgot$ such that
\begin{equation*}
F_o\circ Ad(h)=Ad(h) \circ F_o \text{ for any } h \in H.
\end{equation*}

For this reason, further we will not distinguish an invariant structure $F$ on $G/H$ and its value $F_o$ at the point $o=H$.

\subsection{ Some important classes in generalized Hermitian
geometry}

The concept of generalized Hermitian geometry appeared in the 1980s and is mostly associated with the works of V.F. Kiritchenko (see, for example, \cite{Ki1} and \cite{Ki2}). It should be mentioned that this theory is a natural consequence of the development of Hermitian geometry and the theory of almost contact structures with many applications. 

In the sequel by $\xgot(M)$ we will denote the set of all smooth vector fields on a manifold $M$.

One of the central objects in generalized Hermitian geometry is a {\it metric
$f$-structure} \cite{Ki1} $(f,g)$,  where $f$ is an $f$-structure  compatible with a
(pseudo) Riemannian metric $g=\langle \cdot, \cdot \rangle$ in the
following sense: 
\begin{equation*}\langle fX, Y \rangle + \langle X, fY \rangle
=0 \text{ for any } X,\; Y \in \xgot (M).
\end{equation*}
Evidently, this definition generalizes the notion of an almost Hermitian structure $J$ in Hermitian geometry. A manifold $M$ equipped with a metric $f$-structure is called a {\it metric $f$-manifold}. 

It is worth noticing that the main classes of generalized Hermitian
geometry (see \cite{Ki1}, \cite{B6}, \cite{B2},  \cite{Gr1}, and \cite{Gr2}) in the special case
$f=J$, where $J$ is an almost complex structure, coincide with those of Hermitian geometry (see \cite{GH}).
In this section we will mainly concentrate on the classes ${\bf Kill \, f}$, ${\bf NKf}$, and ${\bf G_1 f}$ of metric $f$-structures.

A fundamental role in generalized Hermitian geometry is played by
the tensor $T$ of type $(2,1)$ which is called a {\it composition
tensor} \cite{Ki1}. In \cite{Ki1} it was shown that such a
tensor exists on any metric $f$-manifold and it is possible to
evaluate it explicitly:
\begin{equation*}
T(X,Y)=\frac{1}{4}f(\nabla _{fX} (f) fY-\nabla _{f^2 X} (f) f^2Y),
\end{equation*}
where $\nabla$ is the Levi-Civita connection of a
(pseudo) Riemannian manifold $(M,g)$, $X,\,Y \in \xgot (M)$.

With the help of this tensor one can define the structure of a so-called {\it adjoint $Q$-algebra} (see \cite{Ki1}) on $\xgot (M)$ by the formula
$X*Y=T(X,Y)$. It gives the opportunity to introduce some classes
of metric $f$-structures in terms of natural properties of the
adjoint $Q$-algebra.

For example, if
\begin{equation}\label{G1f}
T(X,X)=0 \text{ for any } X \in \xgot (M)
\end{equation}
(that is, if $\xgot (M)$ is an anticommutative $Q$-algebra) then $f$ is
referred to as a {\it $G_1 f$-structure}. ${\bf G_1 f}$ denotes
the class of $G_1 f$-structures, which was first introduced (in a more general situation) in \cite{Ki1}.

A metric $f$-structure on $(M,g)$ is said to be a {\it Killing
$f$-structure}  \cite{Gr1, Gr2} if
\begin{equation}\label{Killf}
\nabla _X (f) X=0 \text{ for any } X \in \xgot (M)
\end{equation}
(that is, if $f$ is a Killing tensor). The class of
Killing $f$-structures is denoted by ${\bf Kill\, f}$. 

The defining property of {\it nearly K\"ahler $f$-structures} (or {\it
$NKf$-structures}) is
\begin{equation}\label{NKf}
\nabla _{fX} (f) {fX}=0 \text{ for any } X \in \xgot (M).
\end{equation}
This class of metric $f$-structures, which is denoted by ${\bf
NKf}$, was first determined in \cite{B1} (see also \cite{B6,B3}). It is
not difficult to see that for $f=J$ the classes ${\bf Kill\, f}$ and ${\bf
NKf}$ coincide with the well-known class ${\bf NK}$ of {\it nearly
K\"ahler structures} \cite{G1}.

The following relations between the classes mentioned are evident:
\begin{equation}\label{subsets}
{\bf Kill \, f} \subset  {\bf NKf} \subset {\bf G_1 f}.
\end{equation}

The classical result bellow will be used to rewrite formulas (\ref{G1f}), (\ref{Killf}) and (\ref{NKf}) in a form more suitable for further considerations.

\begin{thm}\cite{KN} Let $(M,g)$ be a Riemannian manifold, $M=G/H$
a homogeneous reductive space with the reductive decomposition
$\ggot=\hgot \oplus \mgot$. Then the Levi-Civita connection with
respect to $g$ can be expressed in the form
\begin{equation}\label{LC}
\nabla_X Y=\frac{1}{2}[X,Y]_{\mgot}+U(X,Y),
\end{equation}
where $U$ is the symmetric bilinear mapping $\mgot \times \mgot
\rightarrow \mgot$ defined by the formula
\begin{equation}\label{U}
2g(U(X,Y),Z)=g(X,[Z,Y]_{\mgot})+g([Z,X]_{\mgot},Y) \text{ for any } X,Y,Z \in
\mgot.
\end{equation}
\end{thm}

It can be shown in the standard way that the application of (\ref{LC}) to
(\ref{G1f}), (\ref{Killf}) and (\ref{NKf}) produces the following
result.
\begin{lem}
Let $(M,g)$ be a Riemannian manifold, $M=G/H$
a reductive homogeneous space with the reductive decomposition
$\ggot=\hgot \oplus \mgot$. Then for an invariant metric $f$-structure $(f,g)$ on $M$ the following holds.
\begin{itemize}

\item[1)] $f \in {\bf G_1f}$ if and only if 
\begin{equation}\label{G1f1}
f(2U(fX,f^2X)-f(U(fX,fX))+f(U(f^2X,f^2X)))=0 \; \text{ for any } \; X \in
\mgot;
\end{equation}

\item[2)] $f \in {\bf NKf}$ if and only if 
\begin{equation}\label{NKf1}
\frac{1}{2}[fX,f^2X]_{\mgot}+U(fX,f^2 X)-f(U(fX,fX))=0 \; \text{ for any }
\; X \in \mgot;
\end{equation}

\item[3)] $f \in {\bf Kill\,f}$  if and only if 
\begin{equation}\label{Killf1}
\frac{1}{2}[X,fX]_{\mgot}+U(X,fX)-f(U(X,X))=0 \; \text{ for any } \; X \in
\mgot.
\end{equation}
\end{itemize}
\end{lem}

\section{Main results}

\begin{as}
Suppose that for a homogeneous reductive space $G/H$ with the reductive decomposition $\ggot=\hgot \oplus \mgot$ the following is true.
\begin{itemize}
\item[$A_1$)] $G$ is a compact semisimple Lie group (hence the Killing form $B$ on $\ggot$ is negative definite).
\item[$A_2$)] The reductive complement $\mgot$ admits the decomposition
\begin{equation}\label{m}
\mgot=\mgot _1 \oplus \mgot _2 \oplus \mgot _3
\end{equation}
into the direct sum of $Ad(H)$-invariant irreducible mutually non-equivalent submodules and
this decomposition is $B$-orthogonal. 
\item[$A_3$)]
\begin{equation}\label{mrel}
0\neq [\mgot _i, \mgot _{i+1}] \subset \mgot _{i+2} \; (\md \;
3),\; i=1,2,3.
\end{equation}
\item[$A_4$)]
\begin{equation}\label{hrel}
[\mgot _i, \mgot _i]\subset \hgot, \;i=1,2,3.
\end{equation}
\end{itemize}
\end{as}

In the view of ${\it A_1)}$ and ${\it A_2)}$ any invariant Riemannian metric $g$ on
$G/H$ is uniquely determined by the triple of positive real
numbers $(a_1,a_2,a_3)$ which implies that
\begin{equation}\label{g}
g=a_1 g_0 \mid _{\mgot _1 \times \mgot _1} + a_2 g_0 \mid _{\mgot _2 \times \mgot _2} + a_3 g_0
\mid _{\mgot _3 \times \mgot _3},
\end{equation}
where $g_0$ is an invariant inner product generated by the negative of the Killing form $B$. Further we will refer to $(a_1,a_2,a_3)$ as to the {\it
characteristic numbers} of $g$. We will also denote the projection
of $X$ onto $\mgot _i$ by $X _i$ for any $X \in \mgot$.

Assumption 1 makes it possible to
calculate the symmetric bilinear mapping $U(X,Y)$ defined in the previous section. The proof of the following result, which is nothing but the simplification of (\ref{U}) in the view of Assumption 1, can be found in \cite{W}.

\begin{lem}
Suppose that $G/H$ satisfies  Assumption 1. Then the symmetric
bilinear mapping $U$ is defined by the formula
\begin{multline}\label{U2}
U(X,Y)=\frac{a_3-a_2}{2a_1}([X_2,Y_3]+[Y_2,X_3])\\
+\frac{a_3-a_1}{2a_2}([X_1,Y_3]+[Y_1,X_3])+\frac{a_2-a_1}{2a_3}([X_1,Y_2]+[Y_1,X_2]).
\end{multline}
\end{lem}

Here and below we assume that $G/H$ satisfies  Assumption 1.
\begin{lem}
For any invariant affinor structure $f$ on $G/H$ $f(\mgot _i)$
$(i=1,2,3)$ is $Ad(H)$-invariant.
\end{lem}
\begin{proof}
${\it A_2)}$ yields that $Ad(h)\mgot _i \subset \mgot _i$ for any
$h \in H$. Hence
\begin{equation*}
f(Ad(H)\mgot _i )\subset f(\mgot _i).
\end{equation*}
$f$ is an invariant affinor structure, therefore
\begin{equation*}
Ad(H)(f(\mgot _i ))\subset f(\mgot _i).
\end{equation*}
\end{proof}
\begin{prop}
Let $f$ be an invariant affinor $f$-structure on $G/H$ with $\Ker f \neq \{0\}$ and $\Image f \neq \{0\}$, $G/H$
satisfies Assumption 1. Then there exists such $i \in \{1,2,3\}$
that either $\Image f=\mgot _i$ or $\Ker f=\mgot _i$.
\end{prop}
\begin{proof}
As $\mgot=\Ker f \oplus \Image
f$, for any $i \in \{1,2,3\}$ we have
\begin{equation*}
\mgot _i=(\Ker f) _i \oplus (\Image f) _i,
\end{equation*}
where
\begin{equation*}
(\Ker f) _i=\mgot_i \cap \Ker f, \; (\Image f) _i=\mgot_i \cap
\Image f.
\end{equation*}

Suppose that there exists $k \in \{1,2,3\}$ such that $(\Ker f) _k \neq
0$ and $(\Image f)_k \neq 0$. Obviously, for any $X \in (\Ker
f)_k$ and $h \in H$
\begin{equation*}
f(Ad(h)X)=Ad(h)(f(X))=0
\end{equation*}
which implies that $(\Ker f)_k$ is $Ad(H)$-invariant.

The same is true for $(\Image f)_k$. Indeed, for any $X \in
(\Image f)_k$ we have $Ad(h)X \in \Image f$ (by Lemma 3) and
$Ad(h)X \in \mgot _k$ (by ${\it A2)}$).

In this way we have obtained that $\mgot _k$ is decomposed into
the sum of the two non-trivial $Ad(H)$-invariant subspaces, which
contradicts Assumption 1. 
\end{proof}
Proposition 1 yields that for any non-trivial invariant affinor $f$-structure $f$ which is not an almost complex structure
 the following is true:
\begin{itemize}
\item [1)] either $f \mid _{\mgot _i}=J$, $f \mid _{\mgot _j \oplus \mgot _k}=0$,
\item [2)] or $f \mid _{\mgot _i}=0 $, $f \mid _{\mgot _j \oplus \mgot _k}=J$,
\end{itemize}
where $\{i,j,k\}= \{1,2,3\}$, $J$ is an almost
complex structure.

Let us consider the first of these two cases. The following statement is valid.
\begin{thm}
Suppose that $G/H$ satisfies Assumption 1, $g$ is an arbitrary
invariant Riemannian metrics on $G/H$. Let $(f,g)$ be an invariant metric $f$-structure and $f \mid _{\mgot
_i}=J$, $f \mid _{\mgot _j \oplus \mgot _k}=0$, where $\{i,j,k\}=\{1,2,3\}$, $J$ is an almost complex structure. Then
\begin{itemize}
\item [1)] $(f,g)$ is not a Killing $f$-structure;
\item [2)] $(f,g)$ belongs to the class ${\bf NKf }$ (and, consequently, to the class ${\bf G_1 f}$).
\end{itemize}
\end{thm}
\begin{proof}
We assume that
\begin{equation}\label{f1}
f \mid _{ \mgot _1}=J, \; f \mid _{\mgot _2 \oplus \mgot _3}=0
\end{equation}
(the results for the other cases are obtained via cyclic
rearrangement of indices). \vspace{5pt}

1) ${\bf Kill \, f}$ is defined by the formula (\ref{Killf1}).
Taking (\ref{U2}), (\ref{f1}) and Assumption 1 into account we obtain
\begin{equation*}
U(X,fX)=\frac{a_3-a_1}{2a_2}[(fX)_1,
X_3]+\frac{a_2-a_1}{2a_3}[(fX)_1,X_2],
\end{equation*}
\begin{equation*}
f(U(X,X))=\frac{a_3-a_2}{a_1}f([X_2,
X_3]).
\end{equation*}
Besides,
\begin{equation*}
\frac{1}{2}[X,fX]_{\mgot}=\frac{1}{2}[X_2,(fX)_1]+\frac{1}{2}[X_3,(fX)_1].
\end{equation*}

Hence, (\ref{Killf1}) is equivalent to the following relation:
\begin{multline*}
\frac{a_3-a_2-a_1}{2a_2}[(fX)_1,
X_3]+\frac{a_2-a_1-a_3}{2a_3}[(fX)_1,
X_2]\\-\frac{a_3-a_2}{a_1}f([X_2, X_3])=0
\end{multline*}
for any  $X \in \mgot$.

By ${\it A3)}$, $[\mgot _i, \mgot _j]\neq 0$ ($i,\,j \in
\{1,2,3\}$, $i \neq j$). Therefore $(f,g)$ belongs to ${\bf Kill\,f}$
if and only if the characteristic numbers of $g$
satisfy the following set of conditions:
\begin{equation*}
\left\{ \begin{aligned}
 \frac{a_3-a_2-a_1}{2a_2}=0, \\
 \frac{a_2-a_1-a_3}{2a_3}=0,\\
 \frac{a_3-a_2}{a_1}=0.
 \end{aligned}
\right.
\end{equation*}
Evidently, this system is inconsistent. \vspace{5pt}

2) The defining property of ${\bf NKf}$ is (\ref{NKf1}). As (\ref{f1}) holds, (\ref{U2}) yields that 
\begin{equation*}
U(fX, fX)=U(fX,f^2X)=0.
\end{equation*}
Moreover, by Assumption 1,
\begin{equation*}
\frac{1}{2}[fX,f^2X]_{\mgot}=\frac{1}{2}[(fX)_1,(f^2X)_1]_{\mgot}=0.
\end{equation*} 

Thus (\ref{NKf1}) holds
for any Riemannian metric. As a particular case, any $f$
satisfying (\ref{f1}) is a $G_1 f$-structure.
\end{proof}

Now let us consider the second group of $f$-structures.

\begin{thm} Suppose that $G/H$ satisfies Assumption 1, $g$ is an arbitrary
invariant Riemannian metrics on $G/H$ with the characteristic numbers $(a_1,a_2,a_3)$. Let $(f,g)$ be an invariant metric $f$-structure, and $f \mid _{\mgot _i}=0$, $f \mid _{\mgot _j \oplus
\mgot _k}=J$, where $\{i,j,k\} = \{1,2,3\}$, $J$ is an
almost complex structure. Then
\begin{itemize}
\item [1)] $(f,g)$ is a $G_1 f$-structure;
\item [2)] $(f,g)$ is a nearly K\"ahler $f$-structure if and only if
$a_j=a_k$ and
\begin{equation}\label{NKfNR}
[fX,f^2X]_{\mgot}=0 \text{ for any } X \in \mgot;
\end{equation}
\item [3)] $(f,g)$ is a Killing $f$-structure if and only if
$a_j=a_k=\frac{4}{3}a_i$ and
\begin{equation*}
\left\{\begin{aligned}
  \,[Z,fZ]_{\mgot}=0,\\
  [Y,fZ]+f([Y,Z])=0
 \end{aligned}
 \right.
\qquad \text{for any }Y \in \mgot _i,\; Z \in \mgot _j \oplus
\mgot _k.
\end{equation*}
\end{itemize}
\end{thm}
\begin{proof}
Without loss of generality it can be assumed that
\begin{equation}\label{f2}
f \mid _{ \mgot _1}=0, \; f \mid _{\mgot _2 \oplus \mgot _3}=J.
\end{equation}

1) It is evident that both
\begin{equation}\label{N1}
U(fX,f^2X)=\frac{a_3-a_2}{2a_1}([(fX)_2,(f^2X)_3]+[(f^2X)_2,(fX)_3])
\end{equation}
and
\begin{equation}\label{N2}
U(fX,fX)=\frac{a_3-a_2}{a_1}[(fX)_2,(fX)_3]
\end{equation}
belong to $\Ker \,f$ for any $X \in \mgot$. Therefore 
(\ref{G1f1}) holds regardless of the choice of $(a_1,a_2,a_3)$.

2) Clearly,
\begin{equation}\label{N3}
\frac{1}{2}[fX,f^2X]_{\mgot}=\frac{1}{2}[(fX)_2,(f^2X)_3]+\frac{1}{2}[(fX)_3,(f^2X)_2].
\end{equation}
Using (\ref{N1}), (\ref{N2}) and (\ref{N3}) we can rewrite
(\ref{NKf1}) as follows:
\begin{equation*}
\frac{a_3-a_2+a_1}{2a_1}[(fX)_2,(f^2X)_3]+\frac{a_3-a_2-a_1}{2a_1}[(f^2X)_2,(fX)_3]=0
\; \text{ for any } X \in \mgot.
\end{equation*}
$a_3-a_2+a_1 \neq 0$ (otherwise $a_3-a_2-a_1=0$ and hence $a_1=0$). Thus $f \in{\bf NKf}$ with
respect to $(a_1,a_2,a_3)$ if and only if
\begin{equation*}
\left\{
\begin{aligned}
 \,[(fX)_2,(f^2X)_3]=\frac{a_3-a_2-a_1}{a_3-a_2+a_1}[(fX)_3,(f^2X)_2], \\
 \left[(fX)_3,(f^2X)_2\right]=\frac{a_3-a_2-a_1}{a_3-a_2+a_1}[(fX)_2,(f^2X)_3]
 \end{aligned}\right.
 \qquad \text{ for any } X \in \mgot
\end{equation*}
(to obtain the second equation we substitute $X$ for $fX$ in the
first one). It follows in the standard way that
\begin{equation*}
\left\{
\begin{aligned}
 \, \left(1-\left(\frac{a_3-a_2-a_1}{a_3-a_2+a_1}\right)^2\right)[(fX)_2,(f^2X)_3]=0, \\
 \left[(fX)_2,(f^2X)_3\right]=\frac{a_3-a_2-a_1}{a_3-a_2+a_1}[(fX)_3,(f^2X)_2]
  \end{aligned}\right.
 \qquad \text{ for any } X \in \mgot.
\end{equation*}

The first equation yields that
$\frac{a_3-a_2-a_1}{a_3-a_2+a_1}=\pm 1$. As $a_1$, $a_2$ and $a_3$
are positive numbers, we have $a_2=a_3$. Then
\begin{equation*}
[(fX)_2,(f^2X)_3]+[(fX)_3,(f^2X)_2]=0.
\end{equation*}

In the view of (\ref{f2}) and Assumption 1 this means that
$[fX,f^2 X] _{\mgot}=0$ for any $X$ in $\mgot$. Thus 2) is proved.

3) As (\ref{subsets}) holds, here we consider $f$-structures
satisfying (\ref{NKfNR}) and invariant Riemannian metrics with
characteristic numbers $(a_1,a_2,a_2)$ $(a_1,\;a_2>0)$ only.

As above, we check that
\begin{equation*}
U(X,X)=\frac{a_2-a_1}{a_2}[X_1,X_2+X_3],
\end{equation*}
\begin{equation*}
U(X,fX)=\frac{a_2-a_1}{2a_2}[X_1,(fX)_2+(fX)_3].
\end{equation*}
Since (\ref{NKfNR}) holds,
\begin{multline*}
\frac{1}{2}[X,fX]_{\mgot}=\frac{1}{2}[X_1,(fX)_2+(fX)_3]+\frac{1}{2}[X_2+X_3,(fX)_2+(fX)_3
]_{\mgot}\\=\frac{1}{2}[X_1,(fX)_2+(fX)_3].
\end{multline*}
Thus (\ref{Killf1}) can be represented as follows:
\begin{equation*}
\frac{2a_2-a_1}{2a_2}[X_1,(fX)_2+(fX)_3]-\frac{a_2-a_1}{a_2}f([X_1,X_2+X_3])=0
\;\text{ for any } \; X \in \mgot.
\end{equation*}
For convenience we shall rewrite it in this way:
\begin{equation*}
\frac{2a_2-a_1}{2a_2}[Y,fZ]-\frac{a_2-a_1}{a_2}f([Y,Z])=0
\;\text{ for any } \; Y \in \mgot _1,\; Z \in \mgot _2 \oplus
\mgot _3.
\end{equation*}
Then it follows that
\begin{equation}\label{1}
[Y,fZ]=\frac{2(a_2-a_1)}{2a_2-a_1}f([Y,Z]) \;\text{ for any } \; Y \in
\mgot _1,\; Z \in \mgot _2 \oplus \mgot _3
\end{equation}
($2a_2-a_1 \neq 0$, because otherwise $a_1=a_2=0$). If we replace
$Z$ by $fZ$ in (\ref{1}) and  then apply $f$ to its both sides, we
obtain
\begin{equation}\label{2}
f([Y,Z])=\frac{2(a_2-a_1)}{2a_2-a_1}[Y,fZ] \;\text{ for any } \; Y \in
\mgot _1, \; Z \in \mgot _2 \oplus \mgot _3.
\end{equation}
(\ref{1}) and (\ref{2}) produce the following system of equations
\begin{equation*}
\left\{ \begin{aligned}
 \,\frac{4(a_2-a_1)^2}{2a_2-a_1}=1, \\
 \left[Y,fZ \right]=\frac{2(a_2-a_1)}{2a_2-a_1}f([Y,Z])
 \end{aligned}
 \right.
 \qquad \text{ for any }\; Y \in \mgot _1,\; Z \in \mgot _2 \oplus \mgot _3.
 \end{equation*}
To conclude the proof, it remains to note that this system is
equivalent to
 \begin{equation*}
\left\{ \begin{aligned}
 a_2=\frac{4}{3}a_1,\\
 \left[Y,fZ \right]+f([Y,Z])=0
 \end{aligned}
 \right.
 \qquad \text{ for any }\; Y \in \mgot _1, Z \in \mgot _2 \oplus \mgot _3.
 \end{equation*}
 \end{proof}
 
 \section{Examples}

\subsection{The manifolds of oriented flags} 

In \cite{IJMMS} we considered manifolds of oriented flags
of the form
 \begin{equation}
 SO(n)/ SO(2) \times SO(n-3)\; (n \geq 4)
 \end{equation}
as homogeneous $\Phi$-spaces \cite{BS1} of order $6$.
 We proved that for any $n \geq 4$ the reductive complement
 $\mgot$ of any such space is decomposed into the direct sum $\mgot=\mgot_1 \oplus \mgot_2 \oplus
 \mgot_3$ of irreducible $Ad(H)$-invariant summands. For the canonical $f$-structures
 on this homogeneous $\Phi$-space of order $6$ we obtained the following result (in the notations of \cite{IJMMS}).

\begin{itemize}
\item [1)] For $f_1(\theta)=\frac{1}{\sqrt{3}}(\theta-\theta^5)$
\begin{equation*}
\Image f_1=\mgot_1 \oplus \mgot_2,\; \Ker f_1=\mgot
 _3.
\end{equation*}
\item [2)] For $f_2(\theta)=\frac{1}{2\sqrt{3}}(\theta-\theta^2+\theta^4-\theta^5)$
\begin{equation*}
\Image f_2=\mgot_2 ,\; \Ker f_2=\mgot _1 \oplus \mgot_3.
\end{equation*}
\item [3)] For
$f_3(\theta)=\frac{1}{2\sqrt{3}}(\theta+\theta^2-\theta^4-\theta^5)$
\begin{equation*}
\Image f_3=\mgot_1 ,\; \Ker f_3=\mgot _2 \oplus \mgot_3.
\end{equation*}
\item [4)] For $f_4(\theta)=\frac{1}{\sqrt{3}}(\theta^2-\theta^4)$
\begin{equation*}
\Image f_4=\mgot_1 \oplus \mgot_2,\; \Ker f_4=\mgot
 _3.
\end{equation*}
\end{itemize}

In \cite{IJMMS} it was checked that for any $i \in \{1,2,3,4\}$ $f_i$ is compatible with any invariant Riemannian metric (\ref{g}), where $g_0=-B(X,Y)=-(n-2) \Tr (X \cdot Y)$. 

The application of Theorem 2 immediately gives us that $(f_2,g)$ and $(f_3,g)$ are not
Killing $f$-structures for any invariant Riemannian
metric $g$. Nevertheless, $(f_2,g)$ and $(f_3,g)$ are nearly K\"ahler
$f$-structures (and, hence, $G_1 f$-structures) with respect to
any invariant Riemannian metric $g$.

Taking account of the facts that $[f_1 X,f_1 ^2 X]=0$, $[f_4 X,f_4 ^2 X]\neq 0$, and $[Y,f_1 Z]+f_1([Y,Z])=0$ for
any $X \in \mgot$, $Y \in \mgot _3$, and $Z \in \mgot_1 \oplus \mgot _2$, by Theorem
3, we obtain
\begin{itemize}
\item [1)] $(f_1,g)$ and $(f_4,g)$ are $G_1 f$-structures for any invariant Riemannian metric $g$;
\item [2)] $(f_1,g)$ belongs to ${\bf NKf}$ if and only if the characteristic numbers of
$g$ are $(s,s,t)$ ($s,t>0$); $(f_4,g)$ is not
a nearly K\"ahler $f$-structure for any invariant
Riemannian metric $g$;
\item [3)] $(f_1,g)$ belongs to ${\bf Kill \, f}$ if and only if the characteristic numbers of $g$ are $(3s,3s,4s)$, where $s>0$. $(f_4,g)$ is
not a Killing $f$-structure for any invariant Riemannian metric $g$.
\end{itemize}

The same results where obtained in \cite{IJMMS} by means of direct
calculations.

\subsection{The complex flag manifold} 

All invariant metric $f$-structures on the complex flag manifold $SU(3)/T_{max}$ ($T_{max}$ is a maximal torus of $SU(3)$) were considered in the view of generalized Hermitian geometry in \cite{arxiv1}. Therefore, here we restrict ourselves to mentioning that $SU(3)/T_{max}$ satisfies the conditions of Assumption 1. Hence Theorems 2 and 3 are applicable in this case.

\subsection{The Stiefel manifold}  Let us consider $G/H=SO(4)/SO(2)$ (a Stiefel manifold). Then 
 \begin{equation*}
 \mgot=\left\{\begin{pmatrix} 0 & a & b_1 & b_2 \\
-a & 0 & c_1 & c_2 \\ -b_1 & -c_1  & 0&
0 \\  -b_2 & -c_2 & 0 & 0 \\
\end{pmatrix}:a, b_1, b_2, c_1, c_2 \in \RR \right\}. 
 \end{equation*}
 
It is not difficult to see that the manifold in question satisfies Assumption 1. Indeed, there is a decomposition of $\mgot$ into the sum of three $Ad(H)$-invariant mutually inequivalent irreducible submodules $\mgot=\mgot_1\oplus\mgot_2\oplus\mgot_3$ (see \cite{Ar}), where 
 
\begin{equation*}
\mgot_1=\left\{\begin{pmatrix} 0 & a & 0 & 0 \\
-a & 0 & 0 & 0 \\ 0 & 0 & 0 &
0 \\  0 & 0 & 0 & 0 \\
\end{pmatrix}:a \in \RR \right\},
\end{equation*}

\begin{equation*}
\mgot_2=\left\{\begin{pmatrix} 0 & 0 & b_1 & b_2 \\
0 & 0 & 0 & 0 \\ -b_1 & 0 & 0 &
0 \\  -b_2 & 0 & 0 & 0 \\
\end{pmatrix}:b_1, b_2 \in \RR \right\},
\end{equation*}

\begin{equation*}
\mgot_3=\left\{\begin{pmatrix} 0 & 0 & 0 & 0 \\
0 & 0 & c_1 & c_2 \\ 0 & -c_1 & 0 &
0 \\  0 & -c_2 & 0 & 0 \\
\end{pmatrix}:c_1, c_2\in \RR \right\}.
\end{equation*}

The conditions $A_3)$ and $A_4)$ are easily checked by straightforward calculations. 

Let us consider the following $f$-structures on this manifold:

\begin{equation*}
f_1:\begin{pmatrix} 0 & a & b_1 & b_2 \\
-a & 0 & c_1 & c_2 \\ -b_1 & -c_1  & 0&
0 \\  -b_2 & -c_2 & 0 & 0 \\
\end{pmatrix} \longrightarrow \begin{pmatrix} 0 & 0 & b_2 & -b_1 \\
0 & 0 & 0& 0 \\ -b_2 &  0 & 0&
0 \\  b_1 & 0 & 0 & 0 \\
\end{pmatrix},
\end{equation*}

\begin{equation*}
f_2:\begin{pmatrix} 0 & a & b_1 & b_2 \\
-a & 0 & c_1 & c_2 \\ -b_1 & -c_1  & 0&
0 \\  -b_2 & -c_2 & 0 & 0 \\
\end{pmatrix} \longrightarrow \begin{pmatrix} 0 & 0 & 0 & 0 \\
0 & 0 & c_2& -c1 \\ 0 &  -c_2 & 0&
0 \\  0 & c_1 & 0 & 0 \\
\end{pmatrix},
\end{equation*}

\begin{equation*}
f_3:\begin{pmatrix} 0 & a & b_1 & b_2 \\
-a & 0 & c_1 & c_2 \\ -b_1 & -c_1  & 0&
0 \\  -b_2 & -c_2 & 0 & 0 \\
\end{pmatrix} \longrightarrow \begin{pmatrix} 0 & 0 & b_2 & -b_1 \\
0 & 0 & c_2& -c1 \\ -b_2 &  -c_2 & 0&
0 \\  b_1 & c_1 & 0 & 0 \\
\end{pmatrix},
\end{equation*}

\begin{equation*}
f_4:\begin{pmatrix} 0 & a & b_1 & b_2 \\
-a & 0 & c_1 & c_2 \\ -b_1 & -c_1  & 0&
0 \\  -b_2 & -c_2 & 0 & 0 \\
\end{pmatrix} \longrightarrow \begin{pmatrix} 0 & 0 & b_2 & -b_1 \\
0 & 0 & -c_2& c1 \\ -b_2 &  c_2 & 0&
0 \\  b_1 & -c_1 & 0 & 0 \\
\end{pmatrix}.
\end{equation*}

There is no difficulty in checking that these $f$-structures are invariant and compatible with any invariant Riemannian metric (\ref{g}), where $g_0=-B(X,Y)=-2 \Tr (X \cdot Y)$. 

By Theorem 2, we obtain that both $(f_1,g)$  and $(f_2,g)$, where $g$ is an arbitrary invariant Riemannian metric, belong to the classes $\bf{NKf}$ and $\bf{G_1 f}$, but they are not Killing $f$-structures.

By Theorem 3, we immediately see that $(f_3,g)$ and $(f_4,g)$ are $G_1f$-structures for any invariant Riemannian metric. 

As $f_3$ does not satisfy (\ref{NKfNR}), $(f_3,g)$ in not an $NKf$-structure, and, consequently, not a Killing $f$-structure with respect to any invariant Riemannian metric.

The verification of the respective conditions of Theorem 3 yields that $(f_4,g)$ is an $NKf$-structure if and only if the characteristic numbers of
$g$ are $(s,t,t)$, where $s,t>0$. $(f_4,g)$ belongs to ${\bf Kill \, f}$ if and only if the characteristic numbers of $g$ are $(4s,3s,3s)$, where $s>0$.

\subsection{The quaternionic flag manifold}
To conclude this paper, we consider the example of the quaternionic flag manifold $G/H=Sp(3)/SU(2)\times SU(2) \times SU(2)$, which also satisfies Assumption 1 \cite{W}.
In this case 
\begin{equation*}
\mgot=\left\{\begin{pmatrix} 0 & x & y \\
-\overline{x} & 0 & z \\ -\overline{y} & -\overline{z} & 0 \\ 
\end{pmatrix}:x,y,z \in \HH \right\}, 
 \end{equation*}
 
\begin{equation*}
 \mgot _1=\left\{\begin{pmatrix} 0 & x & 0 \\
-\overline{x} & 0 & 0 \\ 0 & 0 & 0 \\ 
\end{pmatrix}:x \in \HH \right\}, 
 \end{equation*}

\begin{equation*}
\mgot_2=\left\{\begin{pmatrix} 0 & 0 & y \\
0 & 0 & 0 \\ -\overline{y} & 0 & 0 \\ 
\end{pmatrix}:y \in \HH \right\}, 
 \end{equation*}
 
\begin{equation*} 
 \mgot_3=\left\{\begin{pmatrix} 0 & 0 & 0 \\
 0 & 0 & z \\ 0 & -\overline{z} & 0 \\ 
\end{pmatrix}:z \in \HH \right\}.
 \end{equation*}
 
 The following $f$-structures 
 \begin{equation*}
 f|_{\mgot_p}(X)=(a_1 \mathbf i+a_2 \mathbf j+a_3 \mathbf k)X, \; a_1^2+a_2^2+a_3^2=1,\; a_1, a_2, a_3 \in \RR, X \in \mgot _p,
 \end{equation*}
 
\begin{equation*}
 f|_{\mgot_q \oplus \mgot_r}=0, \{p,q,r\}=\{1,2,3\},
 \end{equation*}
are invariant and compatible with any  invariant Riemannian metric (\ref{g}), where $g_0=-\real (B(X,Y))=-8 \real \Tr (X \cdot Y)$, which is checked by direct calculations. Therefore, by Theorem 2, any of these $f$-structures is both $NKf$- and $G_1 f$-structure. At the same time, it is not a Killing $f$-structure with respect to any invariant Riemannian metric.

Also invariant and compatible with any invariant Riemannian metric (\ref{g}) are $f$-structures of the form

\begin{equation}\label{fq1}
f_1: \begin{pmatrix} 0 & x & y \\
-\overline{x} & 0 & z \\ -\overline{y} & -\overline{z} & 0 \\ 
\end{pmatrix} \longrightarrow \begin{pmatrix} 0 & h_1 x & h_2y \\
-\overline{h_1x} & 0 & 0 \\ -\overline{h_2y} & 0 & 0 \\ 
\end{pmatrix}, 
 \end{equation}
where $h_1$, $h_2 \in \HH$ are such that $\mathrm {Re}\, h_1=\mathrm {Re}\, h_2=0, \; |h_1|=|h_2|=1$.

In this case we have 

\begin{equation*}
[f_1X,f_1^2X]_{\mgot}= \begin{pmatrix} 0 & 0 & 0 \\
0 & 0 & \overline{h_1x}y-\overline{x}h_2y\\ 0 & \overline{h_2y}x-\overline{y}h_1 x & 0 \\ 
\end{pmatrix}, 
\end{equation*}
where  
\begin{equation*} X=\begin{pmatrix} 0 & x & y \\
-\overline{x} & 0 & z \\ -\overline{y} & -\overline{z} & 0 \\ 
\end{pmatrix} \in \mgot.
\end{equation*}
For this reason, $[f_1X,f_1^2X]_{\mgot}=0$ for any $X \in \mgot$ if  and only if $h_1=\overline{h_2}=-h_2$. At the same time, there exist such $Y \in \mgot _3$, $Z \in \mgot_1 \oplus \mgot _2$ that, regardless of the choice of $h_1$ and $h_2$, $[Y,fZ]+f([Y,Z])\neq 0$. 

Thus, an invariant metric $f$-structure $(f_1,g)$, where $f_1$ is of the form (\ref{fq1}), $g$ is an arbitrary Riemannian metric, belongs to the class ${\bf G_1 f}$ and does not belong to the class ${\bf Kill\, f}$.  In this case $(f_1,g)$ is an $NKf$-structure if and only if $h_1=-h_2$ and the characteristic numbers of $g$ are $(\lambda, \lambda, \mu)$, where $\lambda, \mu >0$.

Arguing as above, we obtain that for any invariant Riemannian metric $g$ $(f_2,g)$ and $(f_3,g)$ are $G_1 f$ structures and are not $NKf$-structures (and, consequently, not Killing $f$-structures). Here

\begin{equation*}
f_2: \begin{pmatrix} 0 & x & y \\
-\overline{x} & 0 & z \\ -\overline{y} & -\overline{z} & 0 \\ 
\end{pmatrix} \longrightarrow \begin{pmatrix} 0 & h_1 x & 0 \\
-\overline{h_1x} & 0 & h_2z \\ 0 & -\overline{h_2 z} & 0 \\ 
\end{pmatrix}, 
 \end{equation*}
 
\begin{equation*}
f_3: \begin{pmatrix} 0 & x & y \\
-\overline{x} & 0 & z \\ -\overline{y} & -\overline{z} & 0 \\ 
\end{pmatrix} \longrightarrow \begin{pmatrix} 0 & 0 & h_1y \\
0 & 0 & h_2z \\ -\overline{h_1y} & -\overline{h_2z} & 0 \\ 
\end{pmatrix}, 
 \end{equation*}
 where $h_1$, $h_2 \in \HH$ are such that $\mathrm {Re}\, h_1=\mathrm {Re}\, h_2=0, \; |h_1|=|h_2|=1$.

\newpage
\renewcommand{\refname}{REFERENCES }


\begin{thebibliography}{99}\label{Refs}

\bibitem{V}
D. V. Alekseevsky, V. V. Lychagin, and A. M. Vinogradov, {\it Main
concepts and notions of differential geometry}, Itogi Nauki i
Tekhniki: Fundamental Research VINITI {\bf 28}, 5--289 (1988) (in Russian).

\bibitem{Ar} A. Arvanitoyeorgos, {\it An introduction to Lie groups and the
geometry of homogeneous spaces}, Student Mathematical Library, vol. 22, American Mathematical Society,
Rhode Island, 2003.

\bibitem{B1}
V. V. Balashchenko, {\it Riemannian geometry of canonical
structures on regular $\Phi$-spaces}, Preprint No.174/1994,
Fakult\"at f\"ur Mathematik der Ruhr-Univerit\"at Bochum,
1--19, 1994.

\bibitem{B4} V. V. Balashchenko, {\it Naturally reductive
Killing $f$-manifolds}, Uspekhi Mat. Nauk {\bf 54}, (3),
151--152, (1999); {\sffamily English translation}: Russian Math. Surveys {\bf 54}, (3), 623--625, (1999).

\bibitem{B5}
V. V. Balashchenko, {\it Homogeneous Hermitian $f$-manifolds},
Uspekhi Mat. Nauk {\bf 56}, (3), 159--160, (2001); {\sffamily English translation}: 
Russian Math. Surveys {\bf 56}, (3), 575--577, (2001).


\bibitem{B3}
V. V. Balashchenko, {\it Homogeneous nearly K\"ahler
$f$-manifolds}, Doklady Akademii Nauk {\bf 376}, (4), 439--441, (2001); {\sffamily English translation}: Doklady Mathematics {\bf 63}, (1), 56--58, (2001).

\bibitem{B6}
V. V. Balashchenko, {\it Invariant nearly K\"ahler $f$-structures
on homogeneous spaces}, Contemporary Mathematics {\bf 288},
263--267, (2001).

\bibitem{B2}
V. V. Balashchenko, {\it Invariant structures generated by Lie
group automorphisms on homogeneous spaces}, Proceedings of the
Workshop "Contemporary Geometry and Related Topics" (Belgrade,
Yugoslavia, 15-21 May, 2002). Editors: N.Bokan, M.Djoric,
A.T.Fomenko, Z.Rakic, J.Wess. World Scientific, 1--32, 2004.

\bibitem{arxiv1} V. V. Balashchenko, {\it Invariant f-structures in the generalized Hermitian
geometry}, http: //arxiv.org/abs/math.DG/0503533.

\bibitem{IJMMS} V. V. Balashchenko and A. Sakovich, {\it Invariant f-structures on the flag manifolds $SO(n)/SO(2)\times SO(n-3)$}, Int. J. Math. Math. Sci. {\bf 2006}, 1--15, (2006). {\raggedright

}

\bibitem{BS1}
V. V. Balashchenko and N. A. Stepanov, {\it Canonical affinor structures of classical type
on regular $\Phi$-spaces}, Matematicheskii Sbornik {\bf 186}, (11), 3--34, (1995); {\sffamily English translation}: Sbornik: Mathematics {\bf 186}, (11), 1551--1580, (1995).

\bibitem{G1}
A. Gray, {\it Nearly K\"ahler manifolds}, J. Diff. Geom. {\bf 4}, (3), 283--309, (1970).

\bibitem{GH} A. Gray and L. M. Hervella, {\it The
sixteen classes of almost Hermitian manifolds and their linear
invariants}, Ann. Mat. Pura ed Appl. {\bf 123}, (4), 35--58, (1980).

\bibitem{Gr1}
A. S. Gritsans, {\it Geometry of Killing $f$-manifolds}, Uspekhi
Mat. Nauk {\bf 45}, (4), 149--150, (1990); {\sffamily English translation}: Russian
Math. Surveys {\bf 45}, (4), 168--169, (1990).

\bibitem{Gr2}
A. S. Gritsans, {\it On construction of Killing $f$-manifolds},
Izv. Vyssh. Uchebn. Zaved. Mat. {\bf 6}, 49-57, (1992); {\sffamily English translation}: Soviet Math. (Iz. VUZ) {\bf 36}, (6), (1992).

\bibitem{Ki1} V. F. Kirichenko, {\it Methods of generalized Hermitian geometry
in the theory of almost contact manifolds}, Itogi Nauki i
Tekhniki: Probl. Geom. VINITI {\bf 18}, 25--71, (1986);
{\sffamily English translation}: J. Soviet Math {\bf 42}, (5), (1988).

\bibitem{Ki2} V. F. Kirichenko, {\it Generalized quasi-Kaehlerian manifolds and axioms of CR-submanifolds in generalized Hermitian geometry, I}, Geom.Dedicata {\bf 51}, 75--104, (1994).

\bibitem{KN}
S. Koboyashi and K. Nomizu, {\it Foundations of differential
geometry}, V.2, Intersc. Publ. J.Willey\&Sons, New York-London,
1969.

\bibitem{N}
K. Nomizu, {\it Invariant affine connections on homogeneous
spaces}, Amer. J. Math. {\bf 76}, (1), 33--65, (1954).

\bibitem{W} 
N. Wallach, {\it Compact homogeneous Riemannian manifolds with
strictly positive curvature}, Ann. Math. {\bf 96}, 277--295, (1972).

\bibitem{Y}
K. Yano, {\it On a structure defined by a tensor field $f$ of type
$(1,1)$ satisfying $f^3+f=0$}, Tensor {\bf 14}, 99--109, (1963).


\end{thebibliography}
\end{document}